\documentclass[twoside]{article}
\usepackage[dvips]{graphicx}
\usepackage{qic, epsfig}
\usepackage{bm}
\usepackage{amsmath,amssymb}
\usepackage{type1cm}

\renewcommand{\thefootnote}{\fnsymbol{footnote}}  
\DeclareSymbolFont{lettersA}{U}{txmia}{m}{it}
 \DeclareMathSymbol{\Indi}{\mathord}{lettersA}{'211}
 \usepackage{latexsym}
\usepackage{color}
\usepackage{ulem}

\begin{document}
\setlength{\textheight}{8.0truein}    

\normalsize\textlineskip
\thispagestyle{empty}
\setcounter{page}{1}


\vspace*{0.88truein}

\alphfootnote

\fpage{1}
\centerline{\bf Eigenvalues of two-state quantum walks induced by the Hadamard walk}
\vspace*{0.37truein}
\centerline{\footnotesize
Shimpei Endo\footnote{shimpei.endo@nucl.phys.tohoku.ac.jp}}
\vspace*{0.015truein}
\centerline{\footnotesize\it Frontier Research Institute for Interdisciplinary Science, Tohoku University,  \\
   Department of Physics, Faculty of Science, Tohoku University}
\baselineskip=10pt
\centerline{\footnotesize\it   6-3,Aoba, Aramaki-aza, Aobaku, Sendai, Miyagi 980-8578, Japan}
\vspace*{10pt}

\centerline{\footnotesize
Takako Endo\footnote{
\;endo-takako-sr@ynu.ac.jp(e-mail of the corresponding author)}}
\vspace*{0.015truein}
\centerline{\footnotesize\it Department of Applied Mathematics, Faculty of Engineering, Yokohama National University}
\baselineskip=10pt
\centerline{\footnotesize\it 79-5 Tokiwadai, Hodogaya, Yokohama, 240-8501, Japan}
\vspace*{10pt}

\centerline{\footnotesize
Takashi Komatsu\footnote{komatsu@coi.t.u-tokyo.ac.jp}}
\vspace*{0.015truein}
\centerline{\footnotesize\it Department of Bioengineering School of Engineering,
The University of Tokyo}
\baselineskip=10pt
\centerline{\footnotesize\it Bunkyo, Tokyo 113-8656, Japan}
\vspace*{10pt}

\centerline{\footnotesize 
Norio Konno\footnote{konno@ynu.ac.jp}}
\vspace*{0.015truein}
\centerline{\footnotesize\it Department of Applied Mathematics, Faculty of Engineering, Yokohama National University}
\baselineskip=10pt
\centerline{\footnotesize\it 79-5 Tokiwadai, Hodogaya, Yokohama, 240-8501, Japan}
\vspace*{10pt}

\vspace*{0.225truein}
\vspace*{0.21truein}
\begin{abstract}
Existence of the eigenvalues of the discrete-time quantum walks is deeply related to localization of the walks.
We revealed the distributions of the eigenvalues given by the splitted generating function method (the SGF method) of the quantum walks we had treated in our previous studies. 
In particular, we focused on two kinds of the Hadamard walk with one defect models and the two-phase QWs that have phases at the non-diagonal elements of the unitary transition operators. As a result, we clarified the characteristic parameter dependence for the distributions of the eigenvalues with the aid of numerical simulation.

\end{abstract}

\vspace*{1pt}\textlineskip 
\section{Introduction} 

The discrete-time quantum walks (DTQWs) as quantum counterparts of the classical random walks, which play important roles in various fields, have attracted much attention in the past two decades    
 \cite{Kempe2003,Kendon2007,VAndraca2012,Konno2008b,AharonovEtAl2001,AmbainisEtAl2001, ManouchehriWang2013}. 
As the reviews of the DTQWs, the readers may be referred to \cite{Konno2008b,VAndraca2012}, for instance.
One of the characteristic properties of the DTQWs is {\it localization}, which is defined that the probability a walker
is found at a point does not converge to zero even in the long-time limit. It has been known that there are two-state QWs in one dimension that have localization \cite{EndoEtAl2014, EndoEtAl2015, InuiEtAl2004, KonnoEtAl2013}.

Localization of the DTQWs is closely related to the existence of the eigenvalues of the unitary transition operators.
However, there were few results for the study of localization from the viewpoint of the eigenvalues of the unitary operators, though there were many approaches for localization \cite{EndoEtAl2014, EndoEtAl2015,KonnoEtAl2013,IdeEtAl2014,Machida2016}.
As the rare results,  Komatsu and Konno \cite{KomatsuKonno2019} revealed the absolute continuous part of the spectrum of the Hadamard walk which has attracted much attention for a decade \cite{KonnoEtAl2013, Konno2008b, VAndraca2012, Konno2010}, and cleared the area of bounded-type stationary measures, stationary measures with quadratic divergence, and with exponential divergence.

So far, it has not been clarified the type of and the position to insert defects, that influence the distributions of the eigenvalues for the eigenfunctions in $l^2$-space on $\mathbb{Z}$ of our models. Also, the model parameter dependence of the area of the eigenvalues has not been known. In our study, one side of the dependence has become clear with the help of numerical simulation, and the influence of defect and phase parameters. has been revealed.
Our study may help to discuss localization, lead to classify the stationary measures, and construct a relation with the Spectrum scattering theory.
Here one of the significance to study the stationary measures is to clear the correspondence with that of the classical systems, i.e., the Markov chains.

The rest of this paper is organized as follows. The definitions of our DTQWs and the main results are given in Section \ref{main}. The remaining section is devoted to summarize our results.

\section{Definitions of the DTQWs and the main results}
\label{main}
\noindent
In this paper, we consider the DTQWs on $\mathbb{Z}$, where $\mathbb{Z}$ is the set of integers.
 The quantum walker with two coin states $|L\rangle$ and $|R\rangle$ is supposed to locate at each lattice point on $\mathbb{Z}$ by superposition.
 The system is described on a tensor Hilbert space ${\cal H}_{p}\otimes {\cal H}_{C}$. The Hilbert space ${\cal H}_{p}$ alters the positions and
is spanned by the orthogonal normalized basis $\{|x\rangle;x\in\mathbb{Z}\}$. Also, the Hilbert
space ${\cal H}_{C}$ represents the coin states and is spanned by the orthogonal normalized
basis $\{|J\rangle:J=L,R\}$. We are here allowed to define
\[|L\rangle=\begin{bmatrix} 1\\0  \end{bmatrix},\;|R\rangle=\begin{bmatrix} 0\\1 \end{bmatrix},\]
for instance. We call $|L\rangle$ and $|R\rangle$, the left and right chiralities, respectively.

The DTQWs are defined as unitary processes in which
each coin state at each position varies with given unitary operations. The quantum
walker in this paper is also manipulated by unitary operations. 
The system of the DTQW at time $t$ is represented by 
\[|\Psi_{t}\rangle={}^T \![\cdots,|\Psi_{t}(-2)\rangle,|\Psi_{t}(-1)\rangle,|\Psi_{t}(0)\rangle,|\Psi_{t}(1)\rangle,|\Psi_{t}(2)\rangle,\cdots]\in{\cal H}_{p}\otimes {\cal H}_{C},\]
where $|\Psi_t(x)\rangle={}^T \! [\Psi^{L}_t(x),\Psi^{R}_t(x)]$ is the amplitude of the DTQW at time $t$. Here $L$ and $R$ correspond to the left and right chirarities, respectively, and $T$ stands for the transposed operator.

Let us prepare a sequences of $2 \times 2$ unitary matrices ${\cal A}=\{ A_x : x \in \mathbb{Z} \}$ with
\begin{align*}
A_x = 
\begin{bmatrix} 
a_{x} & b_{x} \\ 
c_{x} & d_{x}
\end{bmatrix}.
\end{align*}
Put \[U^{(s)}=S\oplus_{x\in\mathbb{Z}}A_{x},\]
where $S$ is the standard shift operator defined by
\[S=\sum_{x}(|x\rangle\langle x+1|\otimes|L\rangle\langle L|+(|x\rangle\langle x-1|\otimes|R\rangle\langle R|).\]
Then the time evolution is determined by 
\begin{align*}
\Psi_{t+1} (x) = (U^{(s)}\Psi_{t})(x)=P_{x+1} \Psi_t(x+1) + Q_{x-1} \Psi_t (x-1) \quad (x \in \mathbb{Z}),
\end{align*}
where
\begin{align*}
P_x = 
\begin{bmatrix} 
a_{x} & b_{x} \\ 
0 & 0 
\end{bmatrix}, 
\qquad 
Q_x = 
\begin{bmatrix} 
0 & 0 \\ 
c_{x} & d_{x}
\end{bmatrix}
\end{align*}
with $U_x = P_x + Q_x$. Then $P_x$ and $Q_x$ express the left and right movements, respectively. 
\subsection{Model 1: The Wojcik model}

At first, we focus on the Wojcik model, whose unitary transition operators are
\begin{align}\label{wojcik}
\{A_{x}\}_{x\in\mathbb{Z}}=\left\{
\dfrac{1}{\sqrt{2}}\begin{bmatrix} 1&1\\ 1& -1\end{bmatrix}_{ x=\pm1,\pm2,\cdots}, \;
\dfrac{\omega}{\sqrt{2}}\begin{bmatrix} 1&1\\ 1& -1\end{bmatrix}_{x=0}\right\}
\end{align}
with $\omega=e^{2i\pi\phi},$ where $\phi\in(0,1)$.
The Hadamard walk can be given by $\phi\to0$ in Eq.\eqref{wojcik}.
We note that the Wojcik model has a phase $2\pi\phi$ only at the origin.
By using recurrence equations, Wojcik et al. \cite{WojcikEtAl2012} solved the eigenvalue problem.
 Endo et al. \cite{EndoKonno2014} and Endo and Konno \cite{EndoKonno2015} derived the stationary, the time-averaged limit, and the weak limit measures. They discussed localization and weak convergence, respectively.

\begin{figure}[t]
\begin{minipage}{0.5\hsize}
\centerline{\epsfig{file=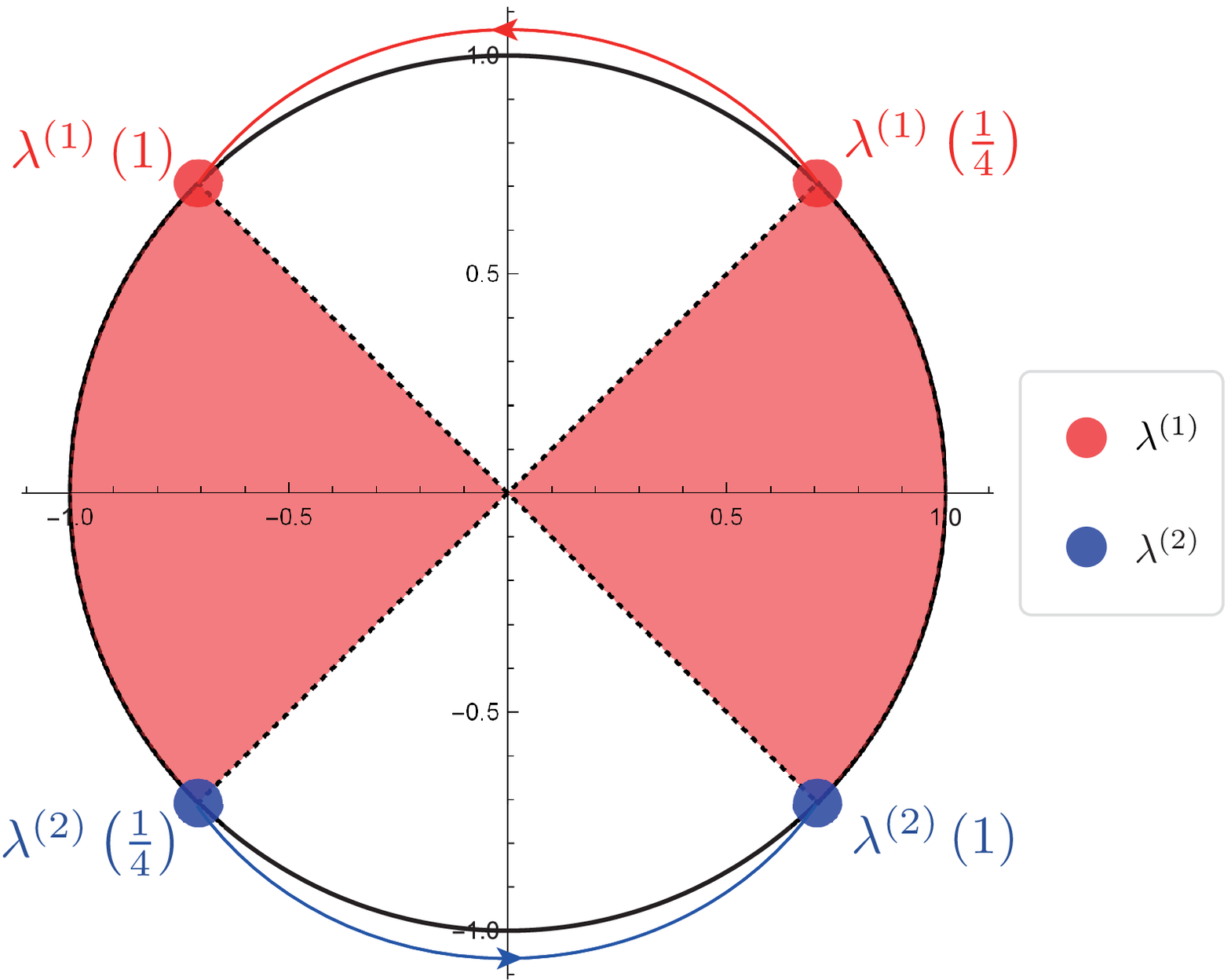, width=7.1cm}}
\vspace*{13pt}
 \end{minipage}
\begin{minipage}{0.5\hsize}
\centerline{\epsfig{file=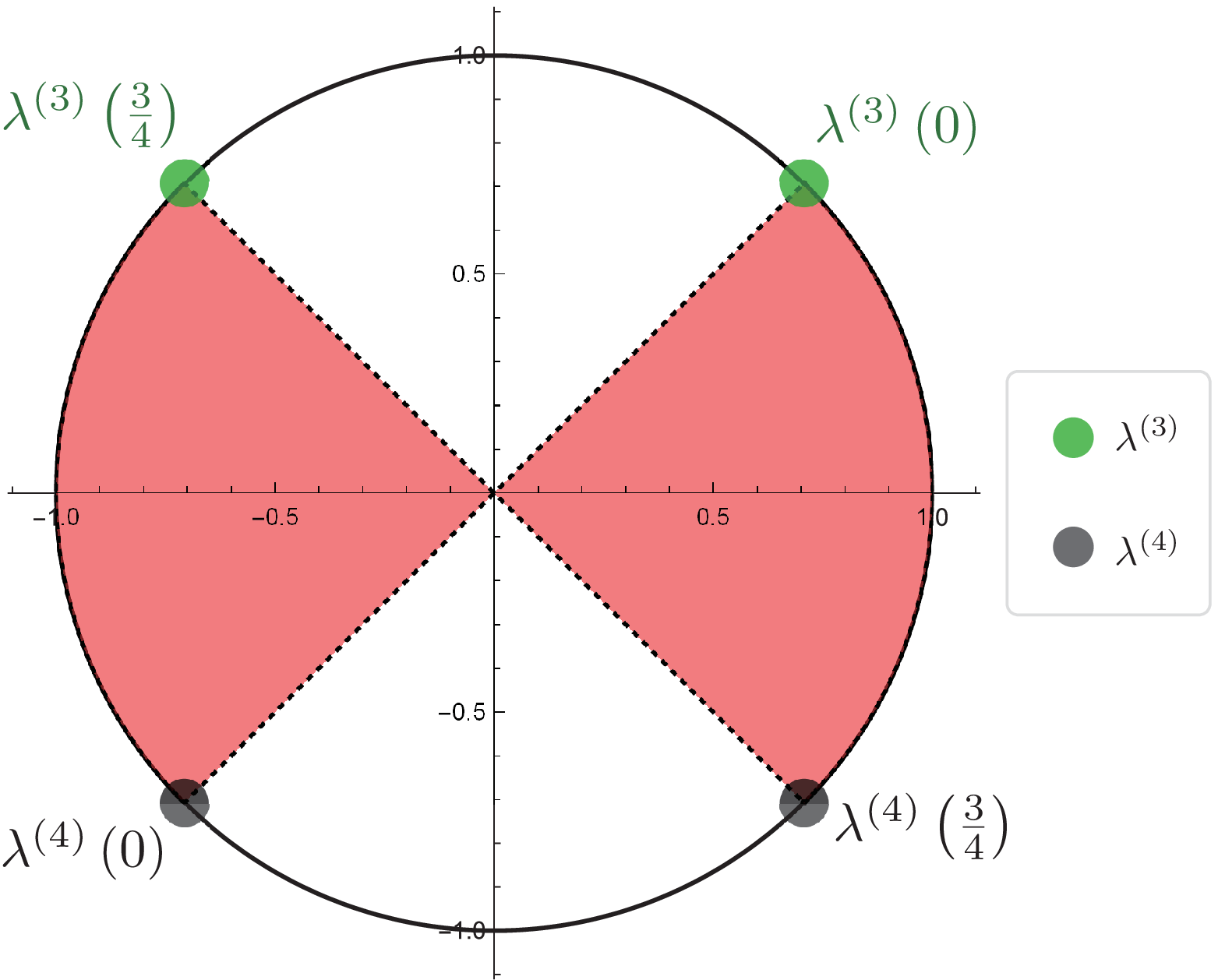, width=7.1cm}} 
\vspace*{13pt}
\end{minipage}\fcaption{The illustrations of the eigenvalues movements of the Wojcik model\\
$(1)\beta=i\alpha$ case. $\lambda^{(1)}(\phi),\lambda^{(2)}(\phi) (\phi\in(\frac{1}{4},1))$\;$(2)\beta=-i\alpha$ case. $\lambda^{(3)}(\phi),\lambda^{(4)}(\phi) (\phi\in(0, \frac{3}{4}))$\\
(The red part is the region of the continuous spectrum of the Hadamard walk, i.e., $\sigma(H)$).}\label{fig.1}
\end{figure}

Now let $\Psi(x)={}^T\![\Psi^{L}(x),\Psi^{R}(x)]$ be the amplitude, and
put $\alpha=\Psi^{L}(0)$ and $\beta=\Psi^{R}(0)$.
Endo et al. \cite{EndoKonno2014} solved the eigenvalue problem
\[U^{(s)}\Psi=\lambda\Psi(\Psi\in Map(\mathbb{Z},\mathbb{C}^{2}), \lambda\in S^{1}),\]
where $S^{1}=\{z\in\mathbb{C};|z|=1\}.$
Here we give the illustrations of the movements of the eigenvalues given by the SGF method \cite{KonnoEtAl2013} of the Wojcik model (Fig.\ref{fig.1}) and a table of the parameter dependence of the eigenvalues  (Table.\ref{table.1}). We remark that each illustration is a diagram of a numerical simulation to investigate the parameter dependence continuously by using mathematica. The SGF method gives the stationary measures with exponential cases and constant case.
Note that the eigenvalues can be obtained by Eqs. (3.8) and (3.9) of Proposition $1$ in \cite{EndoKonno2014}:
\par\indent
\par\noindent
(1) $\beta=i \alpha$ case.
\begin{align*}
\lambda^{2}=\dfrac{\omega(1-2\omega+\omega^{2})-i\omega(1-\omega+\omega^{2})}{1-2\omega+2\omega^{2}}.
\end{align*}
(2) $\beta= -i \alpha$ case.
\begin{align*}
\lambda^{2} 
=\dfrac{\omega(1-2\omega+\omega^{2})+i\omega(1-\omega+\omega^{2})}{1-2\omega+2\omega^{2}}.
\end{align*}
\noindent
Letting 
\[
\left\{ \begin{array}{ll}
\lambda^{(1)}(\phi)=\sqrt{\lambda^{2}},\lambda^{(2)}(\phi)=-\sqrt{\lambda^{2}} & (\beta=i\alpha), \\
\\
\lambda^{(3)}(\phi)=\sqrt{\lambda^{2}},\lambda^{(4)}(\phi)=-\sqrt{\lambda^{2}} & (\beta=-i\alpha), 
\end{array} \right.\]
we specified the regions of the parameter $\phi$ that lead to the eigenfunctions in $l^{2}$-space on $\mathbb{Z}$ by elementary analytic calculations, that is, we have $\phi\in(\frac{1}{4},1)$ for $\beta=i \alpha$ case, and $\phi\in(0, \frac{3}{4})$ for $\beta= -i \alpha$ case.

Now let $\sigma(H)$ be the region of the continuous spectrum of the Hadamard walk. We see that the Wojcik model does not have the eigenvalues on $\sigma(H)$ in the range of the parameter $\phi$.
We notice that despite the two divided cases of the initial state, the distributions of the eigenvalues are the same, and the eigenvalues move allover $S^{1}\setminus\sigma(H)$.

\subsection{Model 2: The Hadamard walk with one defect}
Next, our one-defect model is defined by the set of unitary matrices 
\begin{align}
\{A_{x}\}_{x\in\mathbb{Z}}= 
\left\{ 
\begin{bmatrix}
\cos \xi &  \sin \xi \\
\sin \xi & - \cos \xi 
\end{bmatrix}_{x=0}, \;
H_{x \in \mathbb{Z} \setminus \{0\}}\right\}
\end{align}
with $\xi \in (0, \pi/2)$. We can extend some cases to $\xi = 0$ or $\xi = \pi/2$. Here $H$ is the Hadamard gate defined by
\begin{align*}
H = \dfrac{1}{\sqrt{2}}
\begin{bmatrix}
1 &  1 \\
1 & -1 
\end{bmatrix}.
\end{align*}

\begin{figure}[t]
\begin{minipage}{0.5\hsize}
\centerline{\epsfig{file=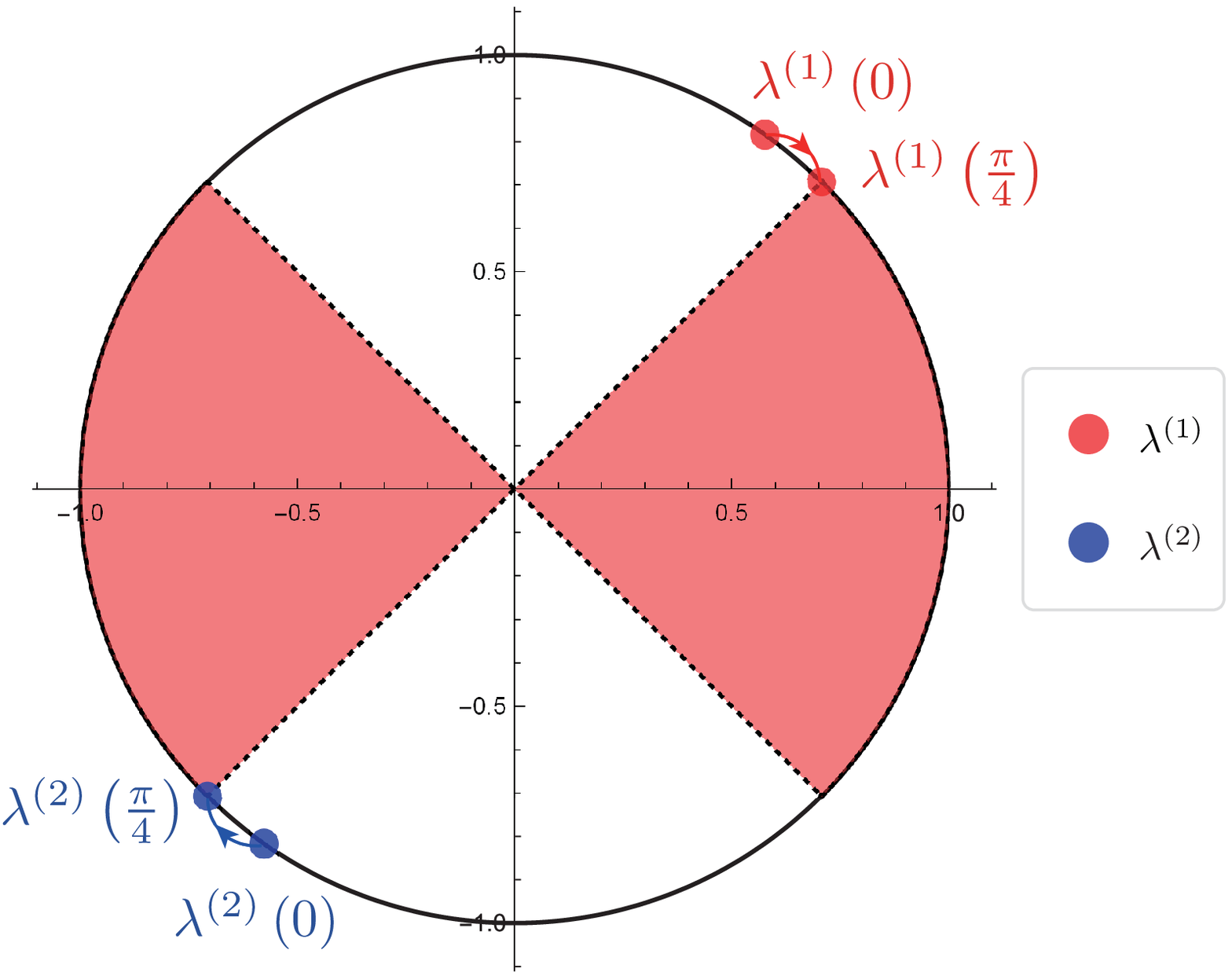, width=7.1cm}}
\vspace*{13pt}
 \end{minipage}
\begin{minipage}{0.5\hsize}
\centerline{\epsfig{file=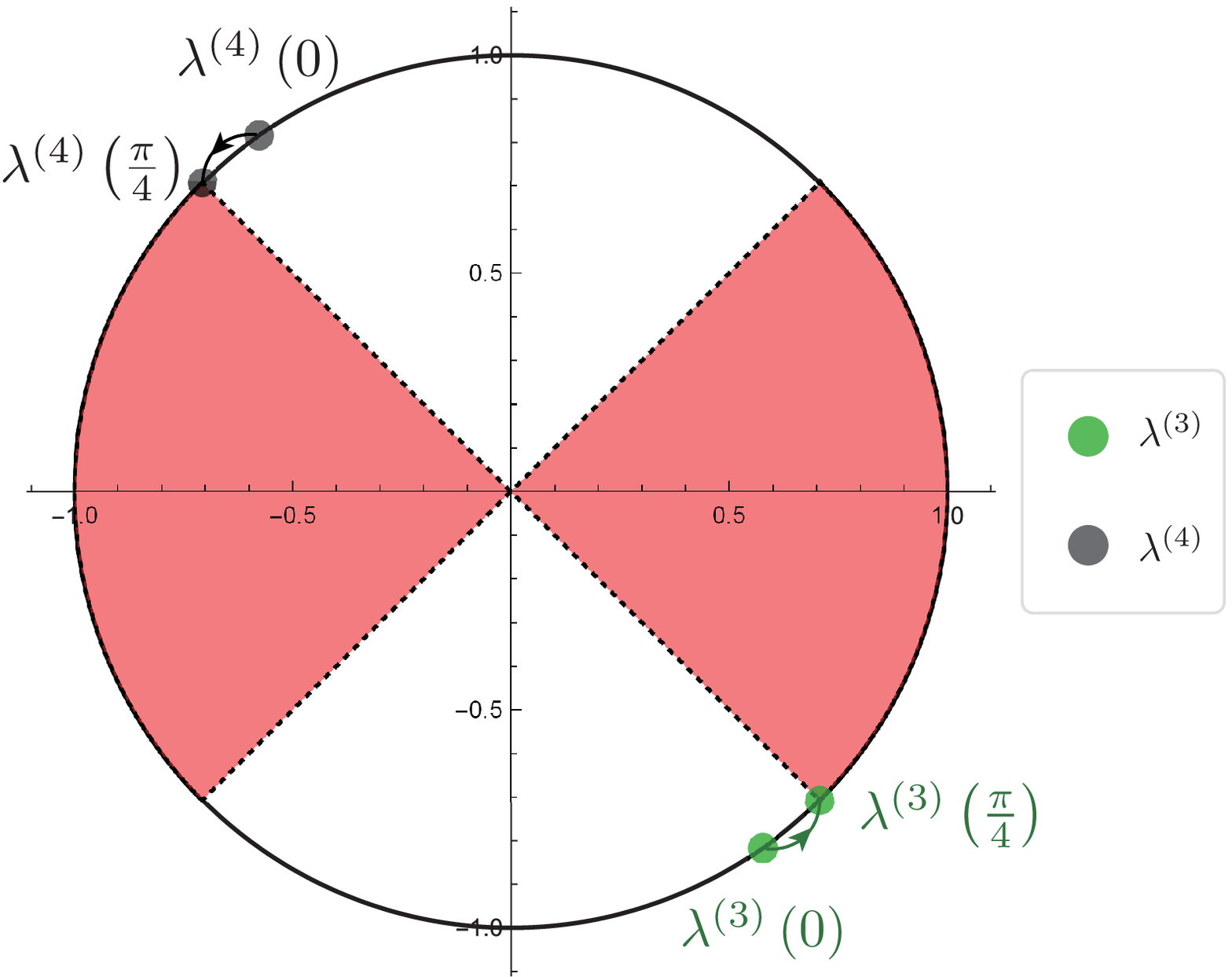, width=7.1cm}} 
\vspace*{13pt}
\end{minipage}\fcaption{The illustrations of the eigenvalues movements of the Hadamard walk with one defect\\
 $(1)\beta=-i\alpha$ case. $\lambda^{(1)}(\xi),\lambda^{(2)}(\xi) (\xi\in(0, \frac{\pi}{4}))$\;$(2)\beta=i\alpha$ case. $\lambda^{(3)}(\xi),\lambda^{(4)}(\xi) (\xi\in(0, \frac{\pi}{4}))$\\
(The red part is the region of the continuous spectrum of the Hadamard walk, i.e., $\sigma(H)$)}\label{fig.2}
\end{figure}

Put $\alpha=\Psi^{L}(0)$ and $\beta=\Psi^{R}(0)$. Then solutions of the eigenvalue problem
\begin{align*}
U^{(s)} \Psi = \lambda \Psi(\Psi\in Map(\mathbb{Z},\mathbb{C}^{2}), \lambda\in S^{1})
\end{align*}
are given in \cite{EndoEtAl2014}. 
Now we show the illustrations of the movements of the eigenvalues given by the SGF method of the Hadamard walk with one defect (Fig.\ref{fig.2}) and a table of the parameter dependence of the eigenvalues (Table.\ref{table.2})
.We remark that each illustration is a diagram of a numerical simulation to investigate the parameter dependence continuously by using mathematica.
Note that the eigenvalues can be obtained by Proposition $3.1$ in \cite{EndoEtAl2014}:
\par\indent
\par\noindent
Put $C=\cos\xi$ and $S=\sin\xi$.\\
(1) $\beta = - i \alpha$ case. We get
\begin{align*}
\lambda = \pm \dfrac{C+(\sqrt{2}-S)i}{\sqrt{3-2 \sqrt{2}S}}.
\end{align*}
\par\noindent
(2) $\beta = i \alpha$ case.  We get
\begin{align*}
\lambda = \pm \dfrac{C-(\sqrt{2}-S)i}{\sqrt{3-2 \sqrt{2}S}}.
\end{align*}
\noindent
Putting
\[
\left\{ \begin{array}{ll}
\lambda^{(1)}(\xi)=\sqrt{\lambda^{2}},\lambda^{(2)}(\xi)=-\sqrt{\lambda^{2}} & (\beta=-i\alpha), \\
\\
\lambda^{(3)}(\xi)=\sqrt{\lambda^{2}},\lambda^{(4)}(\xi)=-\sqrt{\lambda^{2}} & (\beta=i\alpha), 
\end{array} \right.\]
we fined out the regions of the parameter $\xi$ that connect to the eigenfunctions in $l^{2}$-space on $\mathbb{Z}$ by basic analytic calculations, i.e., $\xi\in(o,\frac{\pi}{4})$. 

We emphasize that the eigenvalues emerge only on $S^{1}\setminus\sigma(H)$.
We also notice that the eigenvalues turn in the opposite direction for the two divided cases of the initial state, and the movements of the eigenvalues don't cover $S^{1}\setminus\sigma(H)$.

\subsection{Model 3:The two-phase QW with one defect}
Here we consider the QW whose time-evolution is determined by the unitary transition operators
\begin{align}
\{A_{x}\}_{x\in\mathbb{Z}}=\left\{
\dfrac{1}{\sqrt{2}}\begin{bmatrix}
1 & e^{i\sigma_{+}} \\
e^{-i\sigma_{+}} & -1\\
\end{bmatrix}_{x\geq 1},\;\dfrac{1}{\sqrt{2}}\begin{bmatrix}
1 & e^{i\sigma_{-}} \\
e^{-i\sigma_{-}} & -1\\
\end{bmatrix}_{x\leq -1},\;
\begin{bmatrix}
1 & 0 \\
0 & -1 \\
\end{bmatrix}_{x=0}
\right\}
\label{2-phase.def}\end{align}
with $\sigma_{\pm}\in\mathbb{R}$.
The quantum walker shifts differently in positive and negative parts respectively, and the determinants are independent of the position, that is, $\det(U_{x})=-1$ for $x\in\mathbb{Z}$. The model is called {\it the two-phase QW with one defect} for short.
If $\sigma_{+}=\sigma_{-}$, the model becomes a one-defect QW which has been so far analyzed in detail \cite{KonnoEtAl2013}. We should notice that our model has a defect at the origin, which enables us to analyze the model simply. 

\begin{figure}[h]
\centerline{\epsfig{file=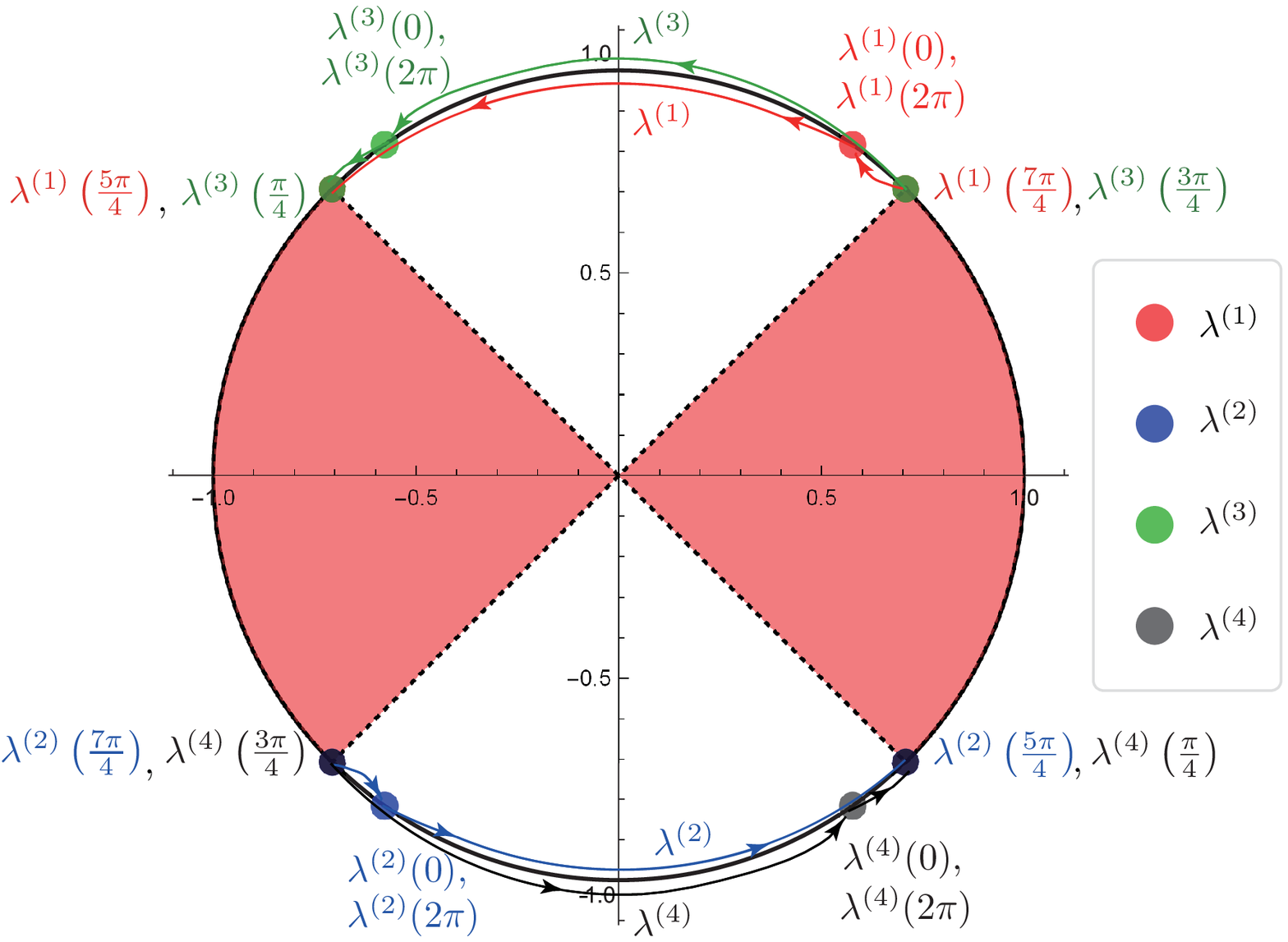, width=14cm}}
\vspace*{13pt}
\fcaption{The illustration of the eigenvalues movements of the two-phase quantum walk with one defect\\
For $\lambda^{(1)}(\sigma)$ and $\lambda^{(2)}(\sigma),$ we have $\sigma\in[0,\frac{5}{4}\pi)\cup(\frac{7}{4}\pi,2\pi].$\\
For $\lambda^{(3)}(\sigma)$ and $\lambda^{(4)}(\sigma),$ we have $\sigma\in[0,\frac{1}{4}\pi)\cup(\frac{3}{4}\pi,2\pi].$\\
(The red part is the region of the continuous spectrum of the Hadamard walk, i.e., $\sigma(H)$)}\label{fig.3}
\end{figure}

Now the solutions of the eigenvalue problem 
\[U^{(s)}\Psi=\lambda\Psi(\Psi\in Map(\mathbb{Z},\mathbb{C}^{2}), \lambda\in S^{1})\]
are described in \cite{EndoEtAl2015}.
Here we put on the illustration of the movements of the eigenvalues given by the SGF method of the two-phase QW with one defect (Fig.\ref{fig.3}) and a table of the parameter dependence of the eigenvalues (Table.\ref{table.3}).
We remark that the illustration is a diagram of a numerical simulation to investigate the parameter dependence continuously by using mathematica.
Note that the eigenvalues can be obtained by Proposition $1$ in \cite{EndoKonno2015}:
\par\indent
\par\noindent
we
\begin{enumerate}
\item $\lambda^{(1)}=\dfrac{\cos\sigma+(\sin\sigma+\sqrt{2})i}{\sqrt{3+2\sqrt{2}\sin\sigma}},\;\lambda^{(2)}=-\lambda^{(1)}.$ 
\item $\lambda^{(3)}=-\dfrac{\cos\sigma+(\sin\sigma-\sqrt{2})i}{\sqrt{3-2\sqrt{2}\sin\sigma}},\;\lambda^{(4)}=-\lambda^{(3)}.$
\end{enumerate}
\noindent
Letting
\[
\lambda^{(1)}(\sigma):=\lambda^{(1)}, \lambda^{(2)}(\sigma):=\lambda^{(2)}, \lambda^{(3)}(\sigma):=\lambda^{(3)}, \lambda^{(4)}(\sigma):=\lambda^{(4)},\]
we specified the regions of the parameter $\sigma$ that lead to the eigenfunctions in $l^{2}$-space on $\mathbb{Z}$ by elementary analytic calculations, that is, we have $\sigma\in[0,\frac{5}{4}\pi)\cup(\frac{7}{4}\pi,2\pi]$ for
$\lambda^{(1)}(\sigma)$ and $\lambda^{(2)}(\sigma),$ and $\sigma\in[0,\frac{1}{4}\pi)\cup(\frac{3}{4}\pi,2\pi]$ for $\lambda^{(3)}(\sigma)$ and $\lambda^{(4)}(\sigma).$

We see that the two-phase quantum walk with one defect does not have the eigenvalues on $\sigma(H)$ in the range of the parameter $\phi$.
We also notice that $\lambda^{1}(\sigma)$ and $\lambda^{3}(\sigma)$, $\lambda^{2}(\sigma)$ and $\lambda^{4}(\sigma)$ turn in the same direction, respectively, and the eigenvalues move allover $S^{1}\setminus\sigma(H)$.

\subsection{Model 4: The complete two-phase QW}
Lastly, we introduce the QW which does not have defects, whose unitary matrices are
\begin{align}
\{A_{x}\}_{x\in\mathbb{Z}}=\left\{
\dfrac{1}{\sqrt{2}}\begin{bmatrix}
1 & e^{i\sigma_{+}} \\
e^{-i\sigma_{+}} & -1 \\
\end{bmatrix}_{x\geq 0},\;
\dfrac{1}{\sqrt{2}}\begin{bmatrix}
1 & e^{i\sigma_{-}} \\
e^{-i\sigma_{-}} & -1\\
\end{bmatrix}_{x\leq -1}\right\}\label{2phase_def}
\end{align}
with $\sigma_{\pm}\in\mathbb{R}$.
The walker steps differently in the spatial regions $x\geq0$ and $x\leq-1$ with the phase parameters $\sigma_{+}$ and $\sigma_{-}$. The QW does not
have defect at the origin, which is in marked contrast to the two-phase QW
with one defect \cite{EndoEtAl2015, EndoEtAl2016}. Hereafter, we call the QW {\it the complete
two-phase QW}.
Putting $\sigma_{+}=\sigma_{-}=0$, the model becomes the Hadamard walk studied in \cite{EndoEtAl2014, Konno2005, KonnoEtAl2013}. 

\begin{figure}[h]
\centerline{\epsfig{file=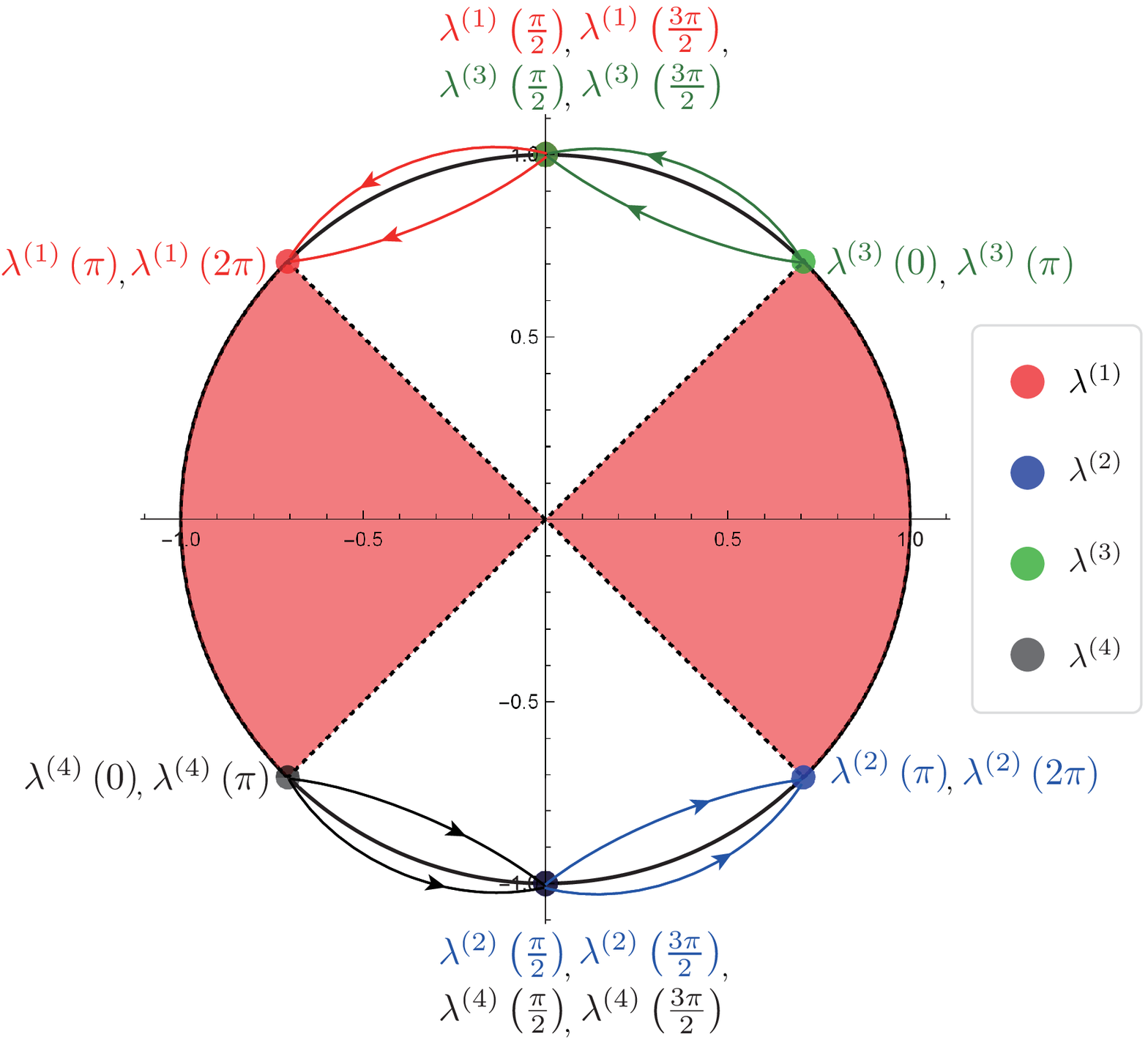, width=14cm}}
\vspace*{13pt}
\fcaption{The illustration of the eigenvalues movements of the complete two-phase quantum walk\\
For $\lambda^{(1)}(\sigma)$ and $\lambda^{(2)}(\sigma),$ we have $\sigma\in[\frac{1}{2}\pi,\pi)\cup(\frac{3}{2}\pi,2\pi].$\\
For $\lambda^{(3)}(\sigma)$ and $\lambda^{(4)}(\sigma),$ we have $\sigma\in[0,\frac{1}{2}\pi)\cup(\pi,\frac{3}{2}\pi].$\\
(The red part is the region of the continuous spectrum of the Hadamard walk, i.e., $\sigma(H)$)}\label{fig.4}
\end{figure}

Let us consider the eigenvalue problem 
\[U^{(s)}\Psi=\lambda\Psi(\Psi\in Map(\mathbb{Z},\mathbb{C}^{2}), \lambda\in S^{1}),\]
whose solutions are given in \cite{EndoEtAl2016}.
Now we show the illustration of the movements of the eigenvalues given by the SGF method of the complete two-phase QW (Fig.\ref{fig.4}) and a table of the parameter dependence of the eigenvalues (Table.\ref{table.4}).
We remark that each illustration is a diagram of a numerical simulation to investigate the parameter dependence continuously by using mathematica.
Note that the eigenvalues can be obtained by Theorem $1$ in \cite{EndoEtAl2016}:

Let $\lambda^{(j)}$ be the eigenvalues of the unitary matrix $U^{(s)}$, and $\Psi^{(j)}(0)$ be the eigenvector at $x=0$, with $j=1,2,3,4$.
Put
\begin{align*}
\left\{
\begin{array}{l}
p=e^{i\sigma_{+}}(e^{-2i\sigma_{-}}-e^{-2i\sigma_{+}}-4e^{-2i\tilde{\sigma}}),\\
q=e^{-2i\sigma_{-}}+e^{-2i\sigma_{+}}-6e^{-2i\tilde{\sigma}},\\
r^{(\pm)}=e^{-i\sigma_{+}}\pm e^{-i\sigma_{-}},
\end{array}
\right
.\end{align*}
where $\tilde{\sigma}=(\sigma_{+}+\sigma_{-})/2$
and
$c\in\mathbb{R}_{+}$ with $\mathbb{R}_{+}=(0, \infty)$. Then we have
\begin{enumerate}
\item $\lambda^{(1)}=\sqrt{\dfrac{p+e^{i\sigma_{+}}r^{(-)}\sqrt{q}}{2(-r^{(-)}-\sqrt{q})}},\;\lambda^{(2)}=-\lambda^{(1)}.$
\item $\lambda^{(3)}=\sqrt{\dfrac{p-e^{i\sigma_{+}}r^{(-)}\sqrt{q}}{2(-r^{(-)}+\sqrt{q})}},\;\lambda^{(4)}=-\lambda^{(3)}.$
\end{enumerate}
\noindent
Putting
\[
\lambda^{(1)}(\sigma):=\lambda^{(1)}, \lambda^{(2)}(\sigma):=\lambda^{(2)}, \lambda^{(3)}(\sigma):=\lambda^{(3)}, \lambda^{(4)}(\sigma):=\lambda^{(4)}, \]
we specified the regions of the parameter $\sigma$ that connect to the eigenfunctions in $l^{2}$-space on $\mathbb{Z}$ by basic analytic calculations, that is, we have $\sigma\in[\frac{1}{2}\pi,\pi)\cup(\frac{3}{2}\pi,2\pi]$ for 
$\lambda^{(1)}(\sigma)$ and $\lambda^{(2)}(\sigma),$ and $\sigma\in[0,\frac{1}{2}\pi)\cup(\pi,\frac{3}{2}\pi]$ for $\lambda^{(3)}(\sigma)$ and $\lambda^{(4)}(\sigma).$

We emphasize that the eigenvalues emerge only on $S^{1}\setminus\sigma(H)$.
We also notice that each eigenvalue turns in a orbit two times, and the movements of the eigenvalues cover allover $S^{1}\setminus\sigma(H)$.

\section{Summary}
In this paper, we focused on four kinds of the QW models originated from the Hadamard walk. As a result, we clarified the characteristic distributions of the eigenvalues given by the SGF method 
Specifically, 
we revealed that the eigenvalues do not emerge on the region of the continuous spectrum of the Hadamard walk.
 Furthermore, by using mathematica, we cleared continuously the parameter dependence of the QW models on the eigenvalues for the eigenfunctions in $l^2$-space on $\mathbb{Z}$. 
We will report on the distributions of the eigenvalues for the three-state QWs we had treated in the forthcoming paper.


\setcounter{footnote}{0}
\renewcommand{\thefootnote}{\alph{footnote}}
\nonumsection{Acknowledgments}
\noindent
S. Endo acknowledges support from  JSPS KAKENHI Grant-in-Aid for Research Activity Start-up 19K21028.
T. Endo is supported by financial support of Postdoctoral Fellowship from Japan Society for the
Promotion of Science.

\par\indent
\par\noindent
\nonumsection{References}

\end{document}